# A Column Generation Algorithm for Vehicle Scheduling and Routing Problems

**Tasnim Ibn Faiz[1], Chrysafis Vogiatzis[2] and Md. Noor-E-Alam[1*]**

[1]Department of Mechanical and Industrial Engineering,
Northeastern University
360 Huntington Ave, Boston, MA 02115, USA

[2]Department of Industrial and Systems Engineering
North Carolina A&T State University
1601 E. Market St, Greensboro, NC, 27411, USA

*Corresponding author email: mnalam@neu.edu

**ABSTRACT**

During natural or anthropogenic disasters, humanitarian organizations face a series of time-sensitive tasks. One of the tasks involves picking up critical resources (e.g., first aid kits, blankets, water) from warehouses and delivering them to the affected people. To successfully deliver these items to the people in need, the organization needs to make decisions that range from the quick acquisition of vehicles from the local market, to the preparation of pickup and delivery schedules and vehicle routes. During crises, the supply of vehicles is often limited, their acquisition cost is steep, and special rental periods are imposed. At the same time, the affected area needs the aid materials as fast as possible, and deliveries must be made within due time. Therefore, it is imperative that the decisions of acquiring, scheduling, and routing of vehicles are made optimally and quickly. In this paper, we consider a variant of a truckload open vehicle routing problem with time windows, which is suitable for modeling vehicle routing operations during a humanitarian crisis. We present two integer linear programming models to formulate the problem. The first one is an arc-based mixed integer linear programming model that is solved using a general purpose solver. The second one, on the other hand, is based on a path-based formulation, for which we design a column generation framework so as to solve it. Finally, we perform numerical experiments and compare the performance of the two models. The comparison shows that the latter path-based formulation outperforms the former without sacrificing solution quality when employing our column generation framework.

**Key words:** Humanitarian logistics, Truckload, Vehicle routing problem, time windows, Column generation, Path generation algorithm.

# 1. Introduction

Our work is largely motivated from the vehicle routing operations that have to take place during a humanitarian crisis from the perspective of a humanitarian logistics organization. In crises arising from natural (flood, earthquakes, hurricanes, etc.) or anthropogenic (e.g., wars) causes, international and local humanitarian organizations (HOs) coordinate towards quickly responding to the needs of the affected area and its population. These organizations are often the only sources available for commodities, such as dry food, water, first aid and other medicine, power sources, clothes. To avoid a duplication of effort, organizations like The Logistics Cluster sometimes act as liaisons for information gathering (data collection), operational decision-making, and coordination of the available resources for improving the timeliness and efficiency of the humanitarian operations.

The number of studies focusing on emergency logistics has seen an increase (Caunhye et al., 2012; Galindo and Batta, 2013) in the last decade after experiencing devastating effects from natural disasters like the Indian Ocean tsunami (2004), Hurricane Katrina in the US (2005), the Haiti earthquake (2010), the earthquake and tsunami in Japan (2011), among others. A number of different assumptions, frameworks, optimization models, and solution approaches have been followed to solve post-disaster logistics problems. We focus on vehicle scheduling and routing for efficient pick-up and drop-off of relief materials from local warehouses to distribution points in a post-disaster setting, which is also the topic of other studies (Lu et al., 2016; Özdamar and Demir, 2012; Özdamar et al., 2004; Wang et al., 2014). Moreover, there have been studies addressing decision problems of relief distribution coupled with evacuation planning (Yi and Özdamar, 2007) and problems of pre-positioning relief materials integrated with post-disaster delivery routing (Caunhye et al., 2016; Rawls and Turnquist, 2012). A common assumption in these studies is that the necessary number of vehicles is readily available for use, which is not often the case in a real-life setting. HOs are not able to maintain a dedicated fleet of vehicles at all locations and immediately dispatch them after a disaster. Instead, when a disaster strikes a specific region, HOs have to make decisions on procuring the necessary transportation means, which include vehicle acquisition with sufficient retention period along with vehicle scheduling and routing. Moreover, after a disaster, vehicle availability becomes limited and their acquisition gets costly. This is why it is crucial for HOs to make prompt and optimal decisions regarding vehicle acquisition, scheduling, and routing to deliver relief materials to the affected people.

To that end, we investigate a special case of pickup and delivery problem with time windows and vehicle rent period restrictions in the context of humanitarian logistics. Our objective is to minimize the overall cost of transportation while ensuring that emergency resources reach the affected areas in due time. The assumptions that we make in our study corresponds to the post-disaster resource delivery premise (see

Section 2.1). We assume that vehicles are tasked to deal with one relief mission at a time and reserve their full capacity for that specific task; this is often the case in humanitarian logistics. In many situations, a vehicle is loaded with multiple products in a single warehouse, then dispatched to the demand location, then sent it to another warehouse for the next mission (Sheu, 2007). Assumption 3 in our study states that each demand location is covered by the nearest depot. This is an assumption followed in (De Angelis et al., 2007) and is common in humanitarian logistics since certain relief items have limited lifetime outside their storage facilities. Finally, assumptions 4 and 5 are common in the literature of general logistics, and also apply in humanitarian operations (see, e.g., the sufficiency assumption in (Berkoune et al., 2012)).The assumptions of our study do not only hold for humanitarian relief operations, but they are also applicable to casualty transportation (see, e.g., (Yi and Özdamar, 2007)), in which case injured people needing medical assistance are transported from their locations to a medical center. This can be easily adapted from our models by appropriately changing the terms of depots and customers to locations with people in need of medical assistance and medical units/hospitals. In the following subsection, we briefly discuss the nature of the problem and the relevant literature.

In land-based transportation, a significant portion of the material volume is moved through truckload deliveries. Each truckload shipment involves a pickup location and a delivery destination, most often with dispatch and arrival time restrictions on the trip ends. In this paper, we study a truckload vehicle scheduling and routing problem with delivery time targets from the perspective of a central resource coordinator. The coordinator has a set of resource depots and is tasked with picking up truckload deliveries from these depots and delivering them to demand locations. The coordinator does not have its own vehicle fleet; rather the vehicles are acquired from a third-party logistics (3PL) service provider when needed. Given a set of orders each with client location with delivery time and its nearest depot location, the supplier faces several challenges, namely i) *vehicle acquisition decisions*: how many vehicles to acquire from dedicated contract carriers and for how long each truck should be rented, ii) *scheduling decisions*: at what time each acquired truck should start its journey and when should it return to the carrier, and iii) *routing decisions*: grouping of orders and assigning a truck to each group. The objective of the supplier is to minimize the total transportation cost through minimizing the number of acquired vehicles and total empty travel distance. The problem is a special case of the vehicle routing problem with time windows (VRPTW) with truckload deliveries and multiple potential starting points for the vehicles routes.

Due to its applicability in improving logistics operations, a number of research studies on VRPTW and its special cases is available in the transportation science literature. Interested readers are referred to the review works in (Braekers et al., 2016; Eksioglu et al., 2009), which provide an excellent overview of the studies in this area. Among the many variants of VRPTW, the pick-up and delivery problem with time windows (PDPTW) (Baldacci et al., 2011), the open vehicle routing problems (OVRP) (Letchford et al., 2007), the dial-a-ride problem (DARP) (Cordeau, 2006), and the vehicle routing problem with crew scheduling (Steinzen et al., 2010) have gained much attention. These studies consider a series of restrictions, and formulate and solve the underlying optimization problems using various techniques. The state-of-the-art for formulating a decision problem for the VRPTW or its special cases is based on arc or path-based set covering models, which are then solved to optimality using branch-and-price or branch-and-cut-and-price using valid inequalities as acceleration techniques. The recent study by Ceselli et al. (Ceselli et al., 2009) can be regarded as representative of current practice in dealing with large-scale VRPTW problems. Their study focuses on the traditional VRPTW with a heterogeneous fleet of vehicles, incompatible products, and driver working hour restrictions. The problem is formulated as a set covering problem and a column generation solution approach is employed to efficiently tackle larger scale instances. The pricing problems, on the other hand, are formulated as resource-constrained shortest path problems, which are solved by a bidirectional dynamic programming algorithm. Ropke and Cordeau (Ropke and Cordeau, 2009) propose a branch-and-cut-and-price algorithm for solving a pickup-and-delivery problem with time windows. They use a set partitioning formulation with two classes of shortest path pricing problems. Several sets of valid inequalities are also incorporated in the branch-and-cut cut framework. A similar problem is studied by Baldacci et al. (Baldacci et al., 2011), in which they devise an exact algorithm based on a set partitioning integer formulation with a bounding procedure. The algorithm is a combination of dual ascent heuristics and a cut-and-column generation procedure that is capable of finding near-optimal solutions. Furtudo et al. (Furtado et al., 2017) study a PDPTW along with pairing and precedence restrictions for pick-up and delivery points. They present a compact formulation and use an exact solution method. A special case of the PDPTW, referred to as pick-up and delivery with transfer for shuttle routes is studied by Masson et al. (Masson et al., 2014). The authors present three models and develop a branch-and-cut-and-price algorithm for solving them. They also present three families of valid inequalities for speeding up the solution process.

Problems involving scheduling of truckload pick-up and delivery, as well as timetabled bus trips, are natural variants of the PDPTW. A vehicle routing and scheduling problem with full truckload (VRPFL) with multiple depot and pick-up time windows is presented by Arunapuram et al. (Arunapuram et al., 2003). The authors formulate the problem as an integer programming model and use a column generation scheme inside a branch-and-bound framework to solve it. Hadjar et al. (Hadjar et al., 2006) use a similar solution approach, to solve a multiple depot vehicle scheduling problem with time windows. In addition to column generation in a branch-and-bound framework, they also use variable fixing, tree search, and cutting planes to speed up the solution process. A multi-vehicle truckload pick-up and delivery problem with time windows (TPDPTW) is studied by Yang et al. (Yang et al., 2004). In this work, the authors consider the costs of empty travel distances of vehicles (travel in which vehicles carry nothing), delay, and rejections of jobs in developing optimal on-line (dynamic) and off-line (stationary) policies. Compared with heuristic rolling-horizon policies, the policy given by re-optimizing the off-line model is found to be more competitive. The authors suggest that taking future job arrival distributions into consideration in scheduling and routing decision making process can result in better policy. A special family of pick-up and delivery problems, e.g., bus trip scheduling, is studied by Kliewer et al. (Kliewer et al., 2006), where a time-space network model is used for formulating the problem. An extension of this work is also available and incorporates crew scheduling in the same modeling framework (Steinzen et al., 2010). As a solution method, a column generation scheme is used, where the pricing subproblems are formulated as constrained shortest path problems in a time-space network. For solving the pricing problems, a dynamic programming algorithm is devised and used.

Among the few studies that consider vehicle tour durations in scheduling and routing decisions is the dial-a-ride problem studied by Cordeau (Cordeau, 2006). This problem is a variant of PDPTW with vehicle capacity and route duration constraints. A branch-and-cut algorithm with several families of valid inequalities is used to solve the model. Archetti and Savelsbergh (Archetti and Savelsbergh, 2009), on the other hand, consider hours of service regulations for drivers in the truckload deliveries with pick-up time windows. A backward search-based algorithm is presented, which can generate a feasible schedule in polynomial time, if one exists. In a more recent study, an exact method for integrated vehicle routing with time windows and truck driver scheduling (VRTDSP) is presented by Goel and Irnich (Goel and Irnich, 2016). The exact branch-and-price algorithm is based on an auxiliary network which enables one to model all possible driver activities on its arcs. This network is then transformed into a parameter-free network which is used in formulating elementary shortest path problems as the pricing problems.

Open vehicle routing problems (OVRP) are a special family of the VRP, where vehicles do not return to a central depot after finishing tours. This set of problems has received increasing attention in recent years due to the rise in 3PL logistics service providers and subscribers. An excellent review of the research works in OVRP is done by Li et al. (Li et al., 2007): the survey provides an overview of the algorithms, datasets, and computational results in this area. One of the first exact methods for OVRP is presented by Letchford et al. (Letchford et al., 2007); they formulate the OVRP as an integer linear programming (ILP) model and develop an algorithm to solve it based on branch-and-cut. Their algorithm is adapted from the branch-and-cut algorithm used for traditional VRP and combined with a separation algorithm involving several sets of valid inequalities. Salari et al. (Salari et al., 2010) present a heuristic procedure for improving an initial feasible solution for an OVRP. This heuristic is based on ILP techniques, where an existing solution is destroyed by randomly removing customers from a path, and then solving an ILP to find an improved solution. Although these studies do not consider pick-up or delivery time windows, they are among the very few studies that attempt to solve the OVRP using an exact method. An OVRP with time windows (OVRPTW) is studied by Repoussis et al. (Repoussis et al., 2007); the authors present a binary ILP model for a generalized OVRP, for which the traditional VRP is a special case. A greedy look-ahead route construction heuristic algorithm is proposed which introduces a new criterion for inserting a delivery location into a vehicle path. This criterion utilizes time windows related information efficiently and keeps track of total waiting time of a vehicle along a partially formed path as well as helps to determine strategies for finding seed customers for building paths. Several other studies investigate the OVRP under various conditions. These studies use heuristic and metaheuristic algorithms to solve the underlying problem, and include hybrid evolution strategy (Repoussis et al., 2010), multi-start algorithms (López-Sánchez et al., 2014), hybrid genetic algorithms (Berger and Barkaoui, 2004), variable neighborhood search algorithms (Fleszar et al., 2009), among others.

In this study, we investigate a special case of the open vehicle pickup and delivery problem from the perspective of a central resource coordinator (the HO) in a post-disaster setting. Following a disaster, the coordinator is tasked with arranging truckload deliveries to the demand locations at fixed delivery times after picking up the load from the nearest resource depots. To meet delivery requirements, the coordinator must rent a sufficient number of vehicles and provide them with tour plans while abiding by all restrictions of vehicle rental periods and delivery times at the client locations. It should be noted here that a vehicle tour neither begins nor ends in a central depot, rather it starts from the first pickup location (resource depot) of the first task on a vehicle tour and ends in the last delivery location (client site) of the last task on the tour. The literature on open vehicle routing problems with time windows (OVRPTW) is scarce; most of these studies adopt heuristic and metaheuristic approaches for solving the underlying optimization problem. The variant that we study here is an open vehicle pickup and delivery problem with fixed delivery times and without a common tour starting point for vehicles. To the best of our knowledge, this particular problem variant has not been investigated fully in the literature, and related problems have been solved using heuristic algorithms. To fill this gap, we present two models to formulate our optimization problems and solve both models to optimality and compare their results. The first model is solved using exact method, while we present fast vehicle path generation algorithms in conjunction with a column generation algorithm to solve the second model.

In Section 2, we propose a mixed integer programming model to formulate the decision problem. Next, in Section 3, we present a column generation scheme to solve the problem that is based on a task conflict graph. We reformulate our problem as a new integer linear program and present two path generation algorithms and a solution algorithm to obtain fast and optimal vehicle acquisition, scheduling, and routing decisions. In Section 4, we compare the computational results obtained from the two approaches. In Section 5, we discuss the limitations of the current study and present potential future extensions.

## 2. Mixed integer linear programming model

The mixed integer linear programming model presented in this section takes the list of orders available at the beginning of the planning horizon as an input. Each order in the list is characterized by the following three features:

1. The client location (destination) and its nearest resource depot (origin).
2. The target delivery time window.
3. The task duration (the sum of loading, unloading, and travelling from origin to destination times).

The model objective is to minimize acquisition, operating, and empty traveling costs of the vehicles, while satisfying delivery targets and vehicle rent period constraints. Our assumptions, decisions, and inputs follow in the next subsections.

### 2.1. Assumptions

The assumptions we make in this work are as follows:

1. A truck is dedicated to perform each task.
2. Loading and unloading times are included in the calculation of task duration.
3. Each client is served from the nearest resource depot.
4. A fraction of the cost of the unused rent period of a divested vehicle is reimbursed by the contract carrier.
5. Supply of resources is sufficient at the depots.

### 2.2. Model decisions

As we will discuss in the mathematical formulation, our decisions are threefold:

1. Number of vehicles to acquire and their lengths of renting periods.
2. Vehicle tours composed of assigned tasks.
3. Starting and ending times of each vehicle tour.

The decisions will be tied to specific decision variables that are introduced later in the section.

## 2.3. Input sets

Our input sets are:

$\mathcal{T}$   Set of tasks indexed by $t$; each task $t$ is characterized by a tuple as follows:

$$\{C_t(client) \quad R_t(nearest\ resource\ depot) \quad T_t(delivery\ time)\}$$

$\mathcal{K}$   Set of vehicles at supplier locations indexed by $k$; each vehicle $k$ is characterized by a tuple as follows:

$$\{r_k(rent\ period) \quad c_k(cost\ per\ unit\ time) \quad p_k(reimbursed\ amount\ per\ unit\ unused\ time)\}$$

To formulate the decision problem as a mixed integer programming (MIP) model, we create a graph $G$ consisting of node set ($\mathcal{N}$) and arc set ($\mathcal{A}$). These sets of nodes and arcs are defined as follows:

1. The node set ($\mathcal{N}$) is obtained by joining the task set ($\mathcal{T}$) with an artificial source node ($D$) and an artificial sink node ($S$).
2. The arc set ($\mathcal{A}$) is composed of three subsets of arcs:
   a. Pull-out arc subset ($A_{po}$): This subset contains pull-out arcs, each of which starts from global source node and ends in a task node, i.e. for all pull-out arc $(i,j) \in A_{po}: i \in D, j \in \mathcal{T}$.
   b. Deadhead arc subset ($A_{dh}$): This subset contains deadhead arcs; a deadhead arc starts from and ends in task nodes, i.e. for all deadhead arc $(i,j) \in A_{dh}: i \in \mathcal{T}, j \in \mathcal{T}, T_i + d_{ij} + \pi_j \leq T_j \leq T_i + d_{ij} + \pi_j + \epsilon$, where $T_i$ and $\pi_i$ are the delivery time and task duration of task $i$ respectively, $d_{ij}$ is the deadhead travel time between task $i \in \mathcal{T}$ and task $j \in \mathcal{T}$ and $\epsilon$ is the maximum allowable vehicle idle time between two consecutive tasks.
   c. Pull-in arc subset ($A_{pi}$): This subset contains pull-in arcs; each of these arcs starts from a task node and ends in global sink node, i.e. for all pull-out arc $(i,j) \in A_{po}: i \in \mathcal{T}, j \in S$.

Figure 1 presents a simplified graphical representation of the problem. The set of arcs with the same color in the figure indicates a representative vehicle tour composed of tasks.

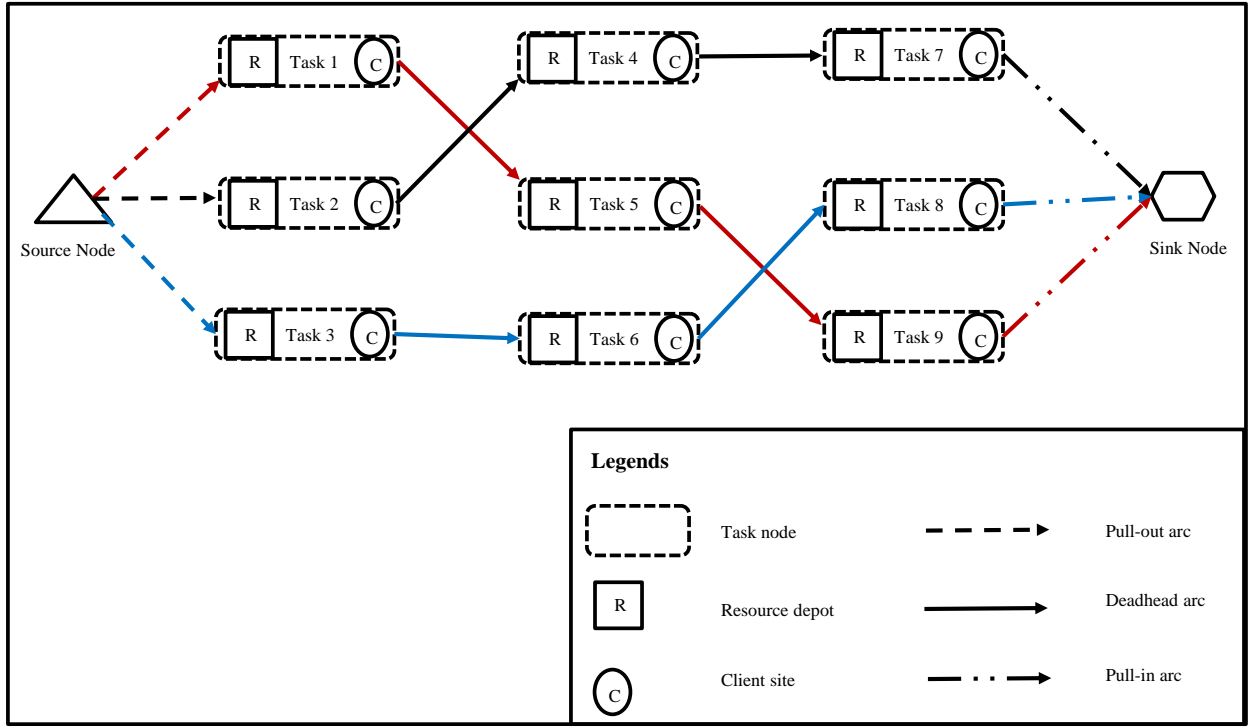

Figure 1: A simplified graphical representation of the network considered in the model.

## 2.4. Mathematical model

Based on the definitions in the previous subsections, our notation can now be presented.

Sets

| | |
|---|---|
| $\mathcal{T}$ | Set of tasks |
| $\mathcal{K}$ | Set of vehicles |
| $\mathcal{A}$ | Set of arcs |

Parameters

| | |
|---|---|
| $r_k$ | Rent period length of vehicle $k \in \mathcal{K}$ |
| $p_k$ | Reimbursed amount for vehicle $k \in \mathcal{K}$ per unit of unused time |
| $c_k$ | Operating cost of vehicle $k \in \mathcal{K}$ per unit of time |
| $d_{ij}$ | Deadhead travel time between task $i \in \mathcal{T}$ and task $j \in \mathcal{T}$ |
| $T_i$ | Delivery time of task $i \in \mathcal{T}$ |
| $\pi_i$ | Duration of task $i \in \mathcal{T}$ |
| $\epsilon$ | Maximum allowable vehicle idle time between two consecutive tasks |

Variables

$Y_{ijk}$     $Y_{ijk} \in \{0,1\}$: $Y_{ijk} = 1$, iff arc $(i,j) \in \mathcal{A}$ is traversed by vehicle $k \in \mathcal{K}$

$\tau_k^s$     Tour starting time for vehicle $k \in \mathcal{K}$

$\tau_k^f$     Tour finishing time for vehicle $k \in \mathcal{K}$

The mathematical model for the decision problem is presented in equations (1)—(13).

$$Minimize \sum_{k \in \mathcal{K}} \sum_{(i,j) \in A_{po}} Y_{ijk}\, r_k c_k (1 - p_k) + \sum_{k \in \mathcal{K}} (\tau_k^f - \tau_k^s) c_k p_k \tag{1}$$

Subject to:

$$\sum_{k \in \mathcal{K}} \sum_{(i,j) \in A_{po} \cup A_{dh}} Y_{ijk} = 1 \qquad \forall j \in \mathcal{T} \tag{2}$$

$$\sum_{(i,j) \in A_{po} \cup A_{dh}} Y_{ijk} - \sum_{(j,i) \in A_{pi} \cup A_{dh}} Y_{jik} = 0 \qquad \forall j \in \mathcal{T}, k \in \mathcal{K} \tag{3}$$

$$\sum_{(i,j) \in A_{po}} Y_{ijk} \leq 1 \qquad \forall k \in \mathcal{K} \tag{4}$$

$$\sum_{(i,j) \in A_{po}} Y_{ijk} - \sum_{(j,i) \in A_{pi}} Y_{jik} = 0 \qquad \forall j \in \mathcal{T}, k \in \mathcal{K} \tag{5}$$

$$\tau_k^s \geq T_j - \pi_j - M(1 - Y_{ijk}) \qquad \forall (i,j) \in A_{po}, k \in \mathcal{K} \tag{6}$$

$$\tau_k^s \leq T_j - \pi_j + M(1 - Y_{ijk}) \qquad \forall (i,j) \in A_{po}, k \in \mathcal{K} \tag{7}$$

$$\tau_k^f \geq T_i - M(1 - Y_{ijk}) \qquad \forall (i,j) \in A_{pi}, k \in \mathcal{K} \tag{8}$$

$$\tau_k^f \leq T_i + M(1 - Y_{ijk}) \qquad \forall (i,j) \in A_{pi}, k \in \mathcal{K} \tag{9}$$

$$T_i \leq T_j + d_{ji} + M(1 - Y_{jik}) \qquad \forall (j,i) \in A_{dh}, k \in \mathcal{K} \tag{10}$$

$$T_i \geq T_j + d_{ji} - M(1 - Y_{jik}) \qquad \forall (j,i) \in A_{dh}, k \in \mathcal{K} \tag{11}$$

$$\tau_k^f - \tau_k^s \leq r_k \sum_{(i,j)\in A_{po}} Y_{jik} \qquad \forall k \in \mathcal{K} \qquad (12)$$

$$\tau_k^s, \tau_k^f \geq 0; \quad Y_{ijk} \in \{0,1\} \qquad \forall (i,j) \in \mathcal{A}, k \in \mathcal{K} \qquad (13)$$

The objective function in (1) minimizes the total transportation cost. If a vehicle $k$ is acquired from a contract carrier with a renting period of $r_k$, then a cost of $r_k c_k$ is incurred. If this vehicle is used for $(\tau_k^f - \tau_k^s)$ units of time, then the supplier receives a reimbursement amount of $r_k - (\tau_k^f - \tau_k^s)p_k$ from the contract carrier. Constraint family (2) ensures that exactly one delivery is made for each task, while constraints (3) restrict that an incoming truck to a task node must depart from that same node. Constraints (4)-(5) ensure that a vehicle tour can only have one starting arc and one ending arc, and they must start from the source node and end in the sink node. The tour start and end times for each vehicle are given by constraints (6)-(9). Feasible deadhead arcs between two consecutive tasks in a vehicle tour are dictated by constraints (10) and (11). A vehicle tour length is constrained by the renting period of that vehicle, which is presented in constraints (12). Finally, the non-negativity and binary nature of the decision variables is restricted with constraints (13). For the parameter $M$ used in constraints (6)-(11), we use the value of the largest delivery time among all the tasks in the task set.

## 3. Column generation algorithm

The MIP model in Section 2 is capable of solving smaller instances (containing a smaller number of tasks) in reasonable time, but it can (expectedly) become computationally expensive as the problem size increases. To solve larger scale problems efficiently using a decomposition technique, we consider a reformulation of the problem that we describe in the subsequent subsections. To facilitate this reformulation, we propose two algorithms for developing potential vehicle paths; each of these paths consists of tasks that can be performed by the same vehicle. To solve this reformulated vehicle routing problem with renting period constraints, column generation is an efficient decomposition approach (Desrosiers et al., 1984). In our approach, first we generate all feasible task combination subsets using Algorithms 1-2 (see subsection 3.1). We also generate sets of potentially good task combinations using Algorithms 3-4 (presented in subsection 3.4). We then design a restricted master problem (presented in subsection 3.2), which we solve considering only the potentially good task combination subsets. If the task combination subsets used in the restricted master problem are not sufficient to ensure optimality, we identify better task combination subsets by solving smaller sub-problems, referred to as the pricing problems (see also subsection 3.3). We now begin our discussion for this approach with algorithms to generate a set of feasible task combination subsets.

### 3.1. Generation of feasible task combination set

We have previously declared $\mathcal{T}$ as the set of tasks and $\mathcal{K}$ as the set of vehicles. Now, we define $C$ as the set of all *task combination subsets* that are *feasible*. To find these subsets, we first construct a graph $\mathcal{G}(\mathcal{V}; \mathcal{A})$ with:

1. Node set $\mathcal{V} = \mathcal{T}$;
2. Arc set $\mathcal{A}$, where an arc $(i, j)$ is in $\mathcal{A}$ if:
   a. $T_j \geq T_i + d_{ij} + \pi_j$ and
   b. $T_j \leq T_i + d_{ij} + \pi_j + \epsilon$.

In this construction, as a reminder, $T_i$ is the delivery time of task $i$, $d_{ij}$ is the deadhead travel time between tasks $i \in \mathcal{T}$ and $j \in \mathcal{T}$, and $\epsilon$ is the maximum allowable time that a vehicle can remain idle between two consecutive tasks. Using graph $\mathcal{G}$, we design Algorithms 1 and 2 to create all feasible task combination subsets (Path) and store them in feasible task set $C$. Algorithm 1 provides the procedure for updating the list of paths (Pathlist) for a task node; this list contains the paths that start from the particular node. This procedure is called upon by Algorithm 2 to perform the procedure for every task node to generate all feasible task combination subsets. These algorithms also provide the duration ($l_{\text{Path}}$) of each task combination subset.

Algorithm 1. $traverse(i, \text{Visited}, \text{Path}, \text{Pathlist})$

1. $(ADJ)_i \leftarrow \{j : j \in \mathcal{T} \text{ and } (i,j) \in \mathcal{A}\}$.
2. Pathlist $\leftarrow$ Pathlist $\cup$ {Path}
3. **for all** $j \in (ADJ)_i$
4.     **if** $j \notin$ Visited **then**
5.         Visited $\leftarrow$ Visited $\cup \{j\}$
6.         Path $\leftarrow$ Path $\cup \{j\}$
7.         $l_{\text{Path}} \leftarrow T_j - T_i + \pi_i$
8.         $traverse(j, \text{Visited}, \text{Path}, \text{Pathlist})$
9.         Visited $\leftarrow$ Visited$\setminus\{j\}$
10.        Path $\leftarrow$ Path$\setminus\{j\}$
11.     **end if**
12. **end for**

Algorithm 2. $feasibleSetgeneration(\mathcal{T},\mathcal{A})$

1. $C \leftarrow \emptyset$
2. **for all** $i \in \mathcal{T}$
3. $\quad Pathlist \leftarrow \emptyset, Path \leftarrow \{i\}, Visited \leftarrow \emptyset$
4. $\quad traverse(i, Visited, Path, Pathlist)$
5. $\quad C \leftarrow C \cup Pathlist$
6. **end for**
7. **return** $C$

Once we have generated set $C$ containing all feasible task combination subsets, we use the following path-based formulation to solve our decision problem.

### 3.2. Restricted master problem (RMP)

Sets

$C$        Feasible task combination set

Parameters

$l_c$        Duration of task combination $c$

$a_{ct}$       1, if task $t$ is in task combination $c$

Variables

$y_{kc}$       1, if vehicle $k$ is assigned to task combination $c \in C$

The restricted master problem is presented below with equations (14)—(18).

$$MP = Minimize \sum_{c \in C} \sum_{k \in \mathcal{K}} y_{kc}[r_k c_k - p_k(r_k - l_c)] \quad (14)$$

Subject to:

$$\sum_{c \in C} \sum_{k \in K} y_{kc} a_{ct} \geq 1 \quad \forall t \in \mathcal{T} \quad (15)$$

$$\sum_{c \in C} y_{kc} \leq 1 \quad \forall k \in \mathcal{K} \quad (16)$$

$$\sum_{c \in C} y_{kc}(l_c - r_k) \leq 0 \quad \forall k \in \mathcal{K} \quad (17)$$

$$y_{kc} \in \{0,1\} \quad \forall k \in \mathcal{K}, c \in C \quad (18)$$

The objective function in (14) minimizes the total cost of vehicle assignment, which is the sum of differences between total vehicle acquisition costs and the total reimbursed amounts received for unutilized portion of rental periods of divested vehicles. Constraints (15) and (16) capture the fact that every vehicle can be assigned to at most one task combination and every task is to be assigned at least to one vehicle. The next constraint family, shown in (17), takes into account each vehicle rental period. This constraint only allows a vehicle to be assigned to a task combination if it can satisfy all the demands in the tasks of that assignment within its rental period. Finally, (18) restricts all of our variables to be binary.

The formulation presented in (14)-(18) has exponentially many variables; hence we cannot hope to directly solve large-scale problems following this formulation within a reasonable amount of time. Instead, we propose the following decomposition scheme.

1. Start with an initial set of task combinations $C_1 \in C$.
2. Define a restricted master problem (RMP) that solves the model in (14)-(18) over set $C_1$ alone.
3. Identify potentially good task combinations to enter in the restricted set $C_1$.
4. Re-solve RMP until we can declare optimality.

The initial restricted set of task combinations $C_1$ will, most probably, not be enough to declare optimality; it is actually possible that it will not even be feasible. Hence, we consider the following restricted dual problem to identify good task combinations to enter in $C_1$:

$$Maximize \sum_{k \in \mathcal{K}} (\lambda_k + \gamma_k) + \sum_{t \in \mathcal{T}} \mu_t \tag{19}$$

Subject to:

$$\lambda_k + \gamma_k + \sum_{t \in T} a_{ct}\mu_t \geq r_k c_k - p_k(r_k - l_c) \qquad \forall k \in \mathcal{K}, c \in C_1, \tag{20}$$

where

- $\mu_t$    are the dual variables for constraints (15)
- $\lambda_k$    are the dual variables for constraints (16)
- $\gamma_k$    are the dual variables for constraints (17).

We are interested in admitting task combinations from the set $C_2 = C \backslash C_1$ that violate constraint (20). Hence, our optimality condition is:

$$\lambda_k + \gamma_k + \sum_{t \in \mathcal{T}} a_{ct}\mu_t - r_k c_k - p_k(r_k - l_c) \geq 0 \qquad \forall k \in \mathcal{K}, c \in C_1 \tag{21}$$

If the above is nonnegative for all vehicles and all task combinations, we can terminate with an optimal solution. For a specific vehicle $\bar{k} \in \mathcal{K}$, we can then solve the pricing problem below to identify a good task combination to feed into the restricted task combination set:

### 3.3. Pricing problem

$$P(\bar{k}): \quad \lambda_k + \gamma_k(l_c - r_k) + Minimize \sum_{c \in C_2} \sum_{t \in \mathcal{T}} a_{ct} \mu_t \phi_c - \sum_{c \in C_2} \phi_c [r_k c_k - p_k(r_k - l_c)] \quad (22)$$

Subject to:

$$\sum_{c \in C_2} \phi_c = 1 \quad (23)$$

$$\phi_c \in \{0,1\} \quad \forall c \in C_2 \quad (24)$$

The pricing sub-problem will give us the most violated constraint (if the objective function is negative) for each vehicle $\bar{k} \in \mathcal{K}$. Iterating over all vehicles, we can get at most $|\mathcal{K}|$ task combinations to enter in our restricted set and re-solve the RMP.

### 3.4. Initial task combination sets (paths) generation for RMP

We mentioned earlier that we cannot solve the RMP in reasonable time if we consider all elements of the task combination set $C$ at once. Instead, we start the RMP with an initial task combination set $C_1$, which is a subset of set $C$. To create the initial task combination set $C_1$, Algorithms 3 and 4 are used, which we describe below. In Algorithm 3, $last(A)$ indicates the last element of set $A$.

Algorithm 3: *greedySelection(i, Taskset)*

1. $S \leftarrow \{j: j \in \mathcal{T} \text{ and } (i,j) \in \mathcal{A}\}$
2. **if** $i = last(\text{Taskset})$ **then do**
3.     **until** $S = \emptyset$
4.         Taskset $\leftarrow$ Taskset $\cup \left\{j: j \in S \text{ and } j = \underset{k}{\operatorname{argmin}}(T_k - T_i + d_{ik} + \pi_k)\right\}$
5.         $S \leftarrow S \setminus last(\text{Taskset})$
6.         $l_{\text{Taskset}} \leftarrow T_{last(\text{Taskset})} - T_{first(\text{Taskset})} + \pi_{first(\text{Taskset})}$
7.     **end**
8. **end if**
9. **return** Taskset, $l_{\text{Taskset}}$

Algorithm 4. *initialSetgeneration($\mathcal{T},\mathcal{A}$)*

1. $C_1 \leftarrow \emptyset$
2. **for all** $i \in \mathcal{T}$ **do**
3.    $S_i \leftarrow \{i\}, l_{S_i} \leftarrow 0$
4.    $(ADJ)_i \leftarrow \{i\} \cup \{j: j \in \mathcal{T} \text{ and } (i,j) \in \mathcal{A}\}$
5.    **if** $(ADJ)_i \neq \emptyset$ **then do**
6.      $S_i \leftarrow greedySelection(i, S_i)$
7.      $C_1 \leftarrow C_1 \cup \{S_i\}$
8.    **end if**
9. **end for**
10. **return** $C_1$

The task combination subsets generated following the above algorithms are stored into the initial task combination set $C_1$. We solve the initial RMP with this restricted task combination set. Rest of the feasible task combinations left in set $C$ (shown in subsection 3.1) are stored in task combination set $C_2$, i.e., $C_2 = C/C_1$. The steps of the solution algorithm under the column generation scheme are defined in Algorithm 5, and in Figure 2, these steps are illustrated schematically.

Algorithm 5: *ColumnGeneration($\mathcal{T},\mathcal{A}$)*

1. $C \leftarrow feasibleSetgeneration(\mathcal{T},\mathcal{A}), C_1 \leftarrow initialSetgeneration(\mathcal{T},\mathcal{A}), C_2 \leftarrow C/C_1$
2. **repeat**
3.   **Solve** RMP without integrality conditions
4.   **Generate** dual variables
5.   **for all** $k \in \mathcal{K}$ **do**
6.     **solve** Pricing problem
7.     **if** $P(k) < 0$ **then do**
8.     $C \leftarrow C_1 \cup \{c: \phi_c = 1\}$
9.     $C_2 \leftarrow C_2 \backslash \{c: \phi_c = 1\}$
10.    **end if**
11.   **end for**
12. **until for all** $\forall k \in \mathcal{K}, P(k) \geq 0$
13. **Solve** RMP with integrality conditions
14. **return** $y_{kc}, MP$

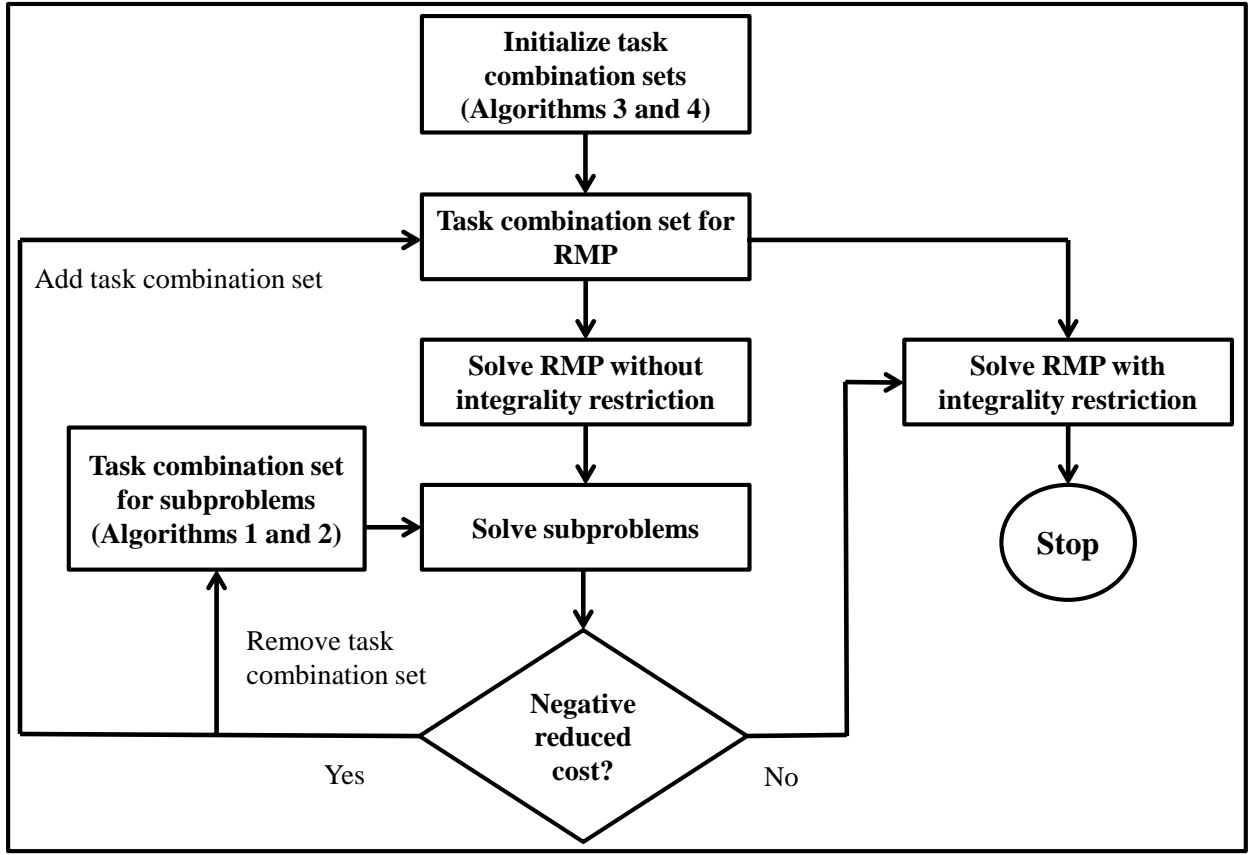

Figure 2: Flow diagram of the column generation algorithm (Algorithm 5).

## 4. Numerical experiments

We implemented the MIP model of (1)—(13) in AMPL (Fourer et al., 1993) and solved it to optimality using CPLEX 12.7.1 solver. We used a computer with Dual Intel Xeon Processor (12 Core, 2.3GHz Turbo) with 64 GB of RAM, and a 64-bit operating system. To implement and solve the master and pricing problems under our column generation scheme, we used the same machine. We then solved the model with randomly generated data for 10, 20, 30, 40, and 50 task problems. Each task has an origin in one of the three resource depot locations and a destination in one of ten potential demand locations. We have considered two classes of vehicles with their respective cost and renting period parameters. The vehicle fleet has the size of 30, where each vehicle in every class is characterized by an operating cost per unit of time, a rental period, and a reimbursement amount per unutilized unit of time for that class. We employed a 'mipgap' setting of 0.001 while solving the model, which indicates that the obtained results are provably within 0.1% of the optimal solution. For each problem size, we solved the model with five sets of randomly generated instances. In every single instance, the column generation algorithm outperforms the exact method when solved to optimality in terms of CPU time, while providing the same (and optimal) objective function values. The summary of our computational experiments is presented in the following Table 1.

Table 1: Comparison between the exact method and column generation algorithm in terms of CPU time.

| Problem size (Number of tasks) | CPU time (in seconds) | | | |
|---|---|---|---|---|
| | Exact method | | Column generation | |
| | Average | Standard deviation | Average | Standard deviation |
| 10 | 1.82 | 1.69 | 1.46 | 0.38 |
| 20 | 4.37 | 3.42 | 2.10 | 0.79 |
| 30 | 44.59 | 47.12 | 4.20 | 2.88 |
| 40 | 251.89 | 118.37 | 9.03 | 3.61 |
| 50 | 910.52 | 350.82 | 17.47 | 4.50 |

## 5. Conclusion

In this study, we investigated a variant of the open vehicle pickup and delivery problem from the perspective of a resource supplier. The supplier is required to rent vehicles from a contract supplier with rental period restrictions to complete a set of pickup and delivery tasks within a target delivery time window. We presented two models for this vehicle acquisition, scheduling, and routing decision problem: an arc-based model and a path-based model. The arc-based MIP model is solved with the exact method; the model provides the number of vehicles to rent, their starting times, and finishing times. The path-based integer model, on the other hand, is solved by a column generation algorithm. Paths (combination of tasks) are generated in the pre-processing steps. We presented two path generation algorithms: the first algorithm generates all feasible paths (compatible task combinations), while the second algorithm provides a set of paths with least idle times for starting the initial restricted master problem. Both algorithms work very fast and the path generation times are very negligible. We compared the results and solution times given by the two models for five different size problems. The path-based model and column generation algorithm outperform the commercial solver. The reduction in solution time by the column generation algorithm is more notable as problem size increases. The reduction in solution time achieved by the column generation algorithm and the fast path generation algorithms indicate that the model can be re-solved repeatedly with updated task list to emulate a dynamic vehicle scheduling and routing problem. During humanitarian crises, when the damage status and demand rates in the affected region are revealed intermittently over time, our algorithm will be most effective in repeatedly solving the decision model with the most updated data to obtain optimal solution very quickly.

As a future extension, we have several potential directions for reformulating our model to make it more applicable in real-world settings and modifying the solution algorithm accordingly. One of these directions includes the consideration of priority in developing paths (task combinations). We have used delivery times as the driving factor in path generation, but it can be easily replaced by other priority rules: for example disaster affected areas with higher percentage of children, older and injured people may get higher priorities, and vehicles may be routed through these areas first. Including uncertainty in i) vehicle availability, ii) vehicle traveling time, and iii) target delivery time in the decision model is another potential extension of our model. To accommodate this uncertainty, we will redesign our path generation algorithms to obtain paths with probabilistic durations and modify our pricing problem to minimize the expected reduced cost of assigning a vehicle to a path. Another possible extension is the integration of vehicle GPS location data with the model; re-solving the model with current vehicle positions and statuses will provide more time efficient and effective decisions.